\documentclass[12pt]{amsart}
\usepackage{amssymb,latexsym}
\usepackage[all]{xy}
\newtheorem{thm}{Theorem}[section]
\newtheorem*{thm*}{Theorem}
\newtheorem{proposition}[thm]{Proposition}
\newtheorem{lemma}[thm]{Lemma}
\newtheorem{cor}[thm]{Corollary}
\theoremstyle{definition}
\newtheorem{remark}[thm]{Remark}
\newtheorem{defn}[thm]{Definition}
 
\newtheorem{notation}[thm]{Notation}  

\newcommand{\by}[1]{\stackrel{#1}{\longrightarrow}}

\newcommand{\rank}{{\rm rank}\,}

\newcommand{\boxtensor}{{\Box\kern-9.03pt\raise1.42pt\hbox{$\times$}}}

\newcommand{\R}{{\mathbb R}}
\newcommand{\Z}{{\mathbb Z}}
\newcommand{\N}{{\mathbb N}}

\newcommand{\tensor}{\otimes}

\newcommand{\sO}{{\mathcal O}}

\renewcommand{\tilde}{\widetilde}
 
\numberwithin{equation}{section}
\newcounter{elno}                

\newcounter{example}[section] 
\def\theexample{\thesection.\arabic{example}}

\begin{document}
\title{Hilbert-Kunz multiplicity and reduction mod~$p$}
\author{V. Trivedi}
\thanks{}
\address{School of Mathematics, Tata Institute of
Fundamental Research,
Homi Bhabha Road, Mumbai-400005, India}
\email{vija@math.tifr.res.in}
\date{}
\maketitle

In this paper, we study the behaviour of  Hilbert-Kunz multiplicities
(abbreviated henceforth to HK multiplicities) of the reductions to
positive characteristics of an irreducible projective curve in
characteristic 0. 

For instance, consider the following question. Let $f$ be a
nonzero irreducible homogeneous element in the polynomial ring
$\Z[X_1,X_2,\ldots,X_r]$, and for any prime number $p\in \Z$, let $R_p =
\Z/p\Z[X_1,X_2, \ldots,X_r]/(f)$ (this is the homogeneous coordinate ring
of a
projective variety over $\Z/p\Z$)). Let $HK(R_p)$ denote the 
Hilbert-Kunz multiplicity of $R_p$ with respect to the graded maximal
ideal. Then one can ask: does $\lim_{p\to \infty}HK(R_p)$ exist? 

 This question was first encountered by the author in a survey article
[C], Problem~4, section~5 (see also Remark~4.10 in [B1]). This seems a
difficult question in
general, as so far, there is no known general formula for HK multiplicity
in terms of `better understood' invariants. There does not seem to even
be a heuristic argument as to why the limit should exist, in general, in
arbitrary dimensions.

However in the case of a projective curve (equivalently 2 dimensional
standard graded ring) over an algebraically closed field of characteristic
$p >0$, one can express HK multiplicity in terms of (i) ``standard'' 
invariants of the curve which are constant in a flat family and
(ii) normalized slopes of the quotients occuring in a strongly semistable
Harder-Narasimhan filtration (HN filtration) of the associated vector
bundle on the curve (see [B1] and [T1]).

Hence, we may pose the question in the following more general
setting. Given a projective curve $X$ defined over a field $k$ of char~$0$
with a vector bundle $V$ on $X$ of rank $r$, there exists a finitely
generated $\Z$-algebra $A$, a projective $A$-scheme $X_A$ such that
  $X_A\tensor_{Q(A)}k = X$, and coherent, locally free sheaves $V_A$ and
$E_{1A},
\ldots, E_{lA}$ on $X_A$ such that,
for all closed points $s \in \mbox{Spec}~A$, if $V_s = V_A\tensor
\overline{k(s)}$, and $E_{i(s)} = E_{iA}\tensor \overline{k(s)}$, then 
 $0\subset E_{1(s)}\subset \cdots
\subset E_{l(s)} \subset V_{s}$ is the
HN filtration of $V_{s}$ (we will
give a detailed version of this in section~2). Choose $s_t\geq 0$ such
that  
$$0\subset F_{1{(s)}} \subset \ldots \subset F_{t_s(s)} \subset
F_{t_s+1(s)} =
F^{s_t*}V_{s}$$
is the strongly semistable HN filtration of $V_s$ (see
Definition~\ref{d5}). 
Denote 
$${\tilde r_{i}}(V_s) =
 \rank(\frac{F_{i(s)}}{F_{i-1(s)}}),~~ 
{\tilde
\mu_{i}}(V_s)
= \mu(\frac{F_{i(s)}}{F_{i-1(s)}}), \mbox{
normalized slope}~~~
{\tilde a_{i}(V_s)} =
\frac{{\tilde
\mu_{i}(V_s)}}{p^{s_t}}.$$
 Let $s_0 \in {\mbox{Spec}~A}$ be the generic
point
of
${\mbox{Spec}~A}$. Then the question is: 

\begin{equation}\label{q2}\mbox{does}~~\lim_{s\to s_0} \sum_i{\tilde
r_{i}}(V_s){\tilde a_i}(V_s)^2 ~~~
\mbox{exist}?
\end{equation}

We  approach the question as follows.
 Following the notation of [L],
for
a vector bundle $V$ on a nonsingular projective curve $X$ in
characteristic $p$, we attach convex polygons as follows. 
Consider the HN filtration 
$$0 =  E_0\subset E_1\cdots \subset E_l\subset E_{l+1} = V.$$
 of $V$. For $k\geq 0$, 
consider the HN filtration 
$$0 = F_0\subset F_1\cdots \subset F_t\subset F_{t+1} = F^{k*}V $$
 of the iterated Frobenius pull back  bundle $F^{k*}V$. 
Let $P(F_i) = (\rank F_i, \deg F_i/p^k)$
in $\R^2$. Let $HNP_{p^k}(V)$ be
the convex
polygon in $\R^2$
obtained by connecting $P(F_0)$, $\ldots$, $P(F_{t+1})$ successively by
line
segments, and connecting the last one with the first one. 
 
Let $p \geq 4(\mbox{genus}(X)-1)(\rank~V)^3$.  Then 
 we prove (Lemma~\ref{l1}) that
the vertices of $HNP_{p^{k-1}}(V)$ are retained as
 a subset of the vertices of $HNP_{p^k}(V)$ and hence $HNP_{p^k}(V) \supset
HNP(V)$. In particular, for $k>>0$, the HN filtration of the bundle
$F^{k*}(V)$ is strongly semistable, therefore 
 Theorem~4.5 of [L] comes as a corollary, in this case.

Now, for every vector
bundle $F_j$ of the HN filtration of $F^{k*}(V)$, if we denote
the slope of the line segment, joining $P(F_{j-1})$ and
$P(F_j)$, by ${\underline \mu}(F_j)/p^k$, and if $E_i$ denotes
 the unique vector bundle occuring in the HN filtration of
$V$ such that $F_j$ `almost descends
to' $E_i$ (see Definition~\ref{d4}), then we prove
(Lemma~\ref{l5}) that 
$${\underline \mu}(F_j)/p^k = 
{\underline \mu}(E_i) + O\left(\frac{1}{p}\right).$$ 
Hence $\lim_{p\to \infty}
{\mbox {Area}}~HNP_{p^k}(V) =  {\mbox {Area}}~HNP(V)$.
In both Lemmas~\ref{l1} and \ref{l5} we make crucial use of 
a result from the paper~[SB] of Shepherd-Barron. 

Now, following the notation set up for the
question~\ref{q2}, we get 
 $$ a_j(V_s) = \frac{{\underline \mu}(F_{j(s)})}{p^s} =
{\underline \mu}(E_{i(s)}) + O\left(\frac{1}{p}\right),$$
where $p = $ characteristic $k(s)$. 
From this we conclude (Proposition~\ref{p1})
that the question~\ref{q2} has an {\it affirmative} answer.

In particular the  Hilbert-Kunz
multiplicities of the
reductions to positive characteristics of an irreducible projective curve
in characteristic 0
have a well-defined limit as the characteristic tends
to $\infty$. This limit,
which is  (relatively) an easier invariant to compute,
  is a {\em lower bound} for the HK multiplicities of the reductions (mod
$p$), though examples of Monsky show that the convergence is {\underline
{not}} monotonic as $p\to \infty$, in general (see Remark~\ref{r7}). 

\section{the HN slope of $F^*V$ in terms of the HN slope 
of $V$}
 Let $X$ be  a nonsingular projective curve of genus $g \geq 1$, over an
algebraically closed
field $k$ of characteristic $p>0$. We recall the following definitions. 

\begin{defn}\label{d1}\begin{enumerate}
\item Let $V$ be a vector bundle on
$X$. We
say $V$ is a {\em semistable} vector bundle on $X$ if, for every
subbbundle
$F \subseteq V$, we have
$$\mu(F) := \frac{\deg~F}{\rank~F} \leq \mu(V).$$

Moreover 
\item $V$ is {\em strongly semistable} if $F^{s*}(V)$ is
semistable for
every $s^{th}$
iterated power of the absolute Frobenius map $F:X\to X$.
\end{enumerate}
\end{defn}

\begin{defn}\label{d2}Let $V$ be a vector bundle on $X$.
A filtration of $V$ by vector subbundles 
\begin{equation}\label{e11} 0 = E_0 \subset E_1\subset \cdots \subset E_l
\subset E_{l+1} =V
\end{equation}
is a {\em Harder-Narasimhan filtration} 
if 
\begin{enumerate}
\item the vector bundles $E_1$, $E_2/E_1$, \ldots, $E_{l+1}/E_l$ are all
semistable.
\item $\mu(E_1) > \mu(E_2/E_1) > \ldots > \mu(E_{l+1}/E_l)$.
\end{enumerate}\end{defn}

\begin{remark}\label{r3}For any Harder-Narasimhan
filtration (we would
call it HN filtration from now onwards), denoted 
as in \ref{e11}, the following is true (see [HN], Lemma~1.3.7),  
\begin{enumerate}
\item the filtration is unique for $V$,
\item $\mu(E_1) > \mu(E_2) >\cdots > \mu(E_{l+1}) = \mu(V)$, 
\item $\mu(E_i/E_{i-1}) \geq \mu(V) \geq \mu(E_{i+1}/E_i)$, for some
$1\leq
i\leq k$.
\end{enumerate}\end{remark}

\begin{notation}\label{n1}If 
$$0=E_0 \subset E_1\subset E_2\subset \cdots \subset E_l\subset E_{l+1}=
V$$ 
is the HN filtration for a vector
bundle $V$ on $X$ then we denote
$${\underline \mu(E_i)} =
\mu\left(\frac{E_i}{E_{i-1}}\right), ~~~
\mu_{\rm{max}}(V) = \mu(E_1)~~~\mbox{and}~~\mu_{\rm{min}}(V) =
\mu\left(\frac{V}{E_l}\right).$$
\end{notation}

\begin{lemma}\label{l0}Let $V$ be a vector bundle over $X$ of rank $r$ and
let 
$$0 ={E_0} \subset E_1\subset {E_2} \subset \cdots
\subset {
E_l}\subset { E_{l+1}} = V$$
be the HN filtration of $V$. Let
$\mu_{i}= \mu({
E_{i}}/{E_{i-1}})$. 
Then $$r^3 > \frac{r-1}{\mu_i-\mu_{i+1}}.$$
\end{lemma}
\begin{proof} Let us denote ${\bar r_i} = \rank~{E_i}/{E_{i-1}}$ and 
${\bar d_i} = \deg~{E_i}/{E_{i-1}}$.
Then 
$$\frac{r-1}{\mu_i-\mu_{i+1}} = \frac{r-1}{{\bar d_i}/{\bar r_1}-{\bar
d_{i+1}}/{\bar r_{i+1}}} = \frac{(r-1){\bar r_i}{\bar r_{i+1}}}{{\bar
d_i}{\bar
r_{i+1}}-{\bar d_{i+1}}{\bar r_i}}.$$
 But 
$$\mu_i-\mu_{i+1} > 0
\implies {\bar d_i}{\bar r_{i+1}}-{\bar d_{i+1}}{\bar r_i} >0 \implies
{\bar d_i}{\bar r_{i+1}}-{\bar d_{i+1}}{\bar r_i} \geq 1.$$
 Therefore 
$$\frac{r-1}{\mu_i-\mu_{i+1}}\leq (r-1){\bar r_i}{\bar r_{i+1}}
\leq r^3.$$
This proves the lemma.\end{proof}

Now we prove the following crucial lemma.

\begin{lemma}\label{l1}Let $V$ be a vector bundle on $X$ as in 
Lemma~\ref{l0}. Assume that ${\rm {char}.}k =  
p > 4(g-1)r^3$.
 Then, 
$$F^*{E_1} \subset F^*{E_2} \subset \cdots \subset 
F^*{E_l} \subset F^*V$$
is a subfiltration of the HN filtration of $F^*V$, that is,
if 
$$0\subset {\tilde E_1}\subset \cdots \subset {\tilde E_{l_1+1}} = F^*V$$
 is the HN
filtration
of $F^*V$ then for every $1\leq i \leq l$  there exists $1\leq j_i\leq
l_1$ such that $F^*{ E_i} = {\tilde E_{j_i}}$. \end{lemma}
\begin{proof}For each $0\leq i\leq l+1$, let 
$$F^*{E_i} \subset E_{i1} \subset \cdots \subset E_{it_i} \subset
F^*{E_{i+1}}$$
be a filtration of vector bundles on $X$ such that 
$$0\subset \frac{E_{i1}}{F^*{E_i}}\subset \frac{E_{i2}}{F^*{
E_i}}\subset
\cdots \subset \frac{F^*{E_{i+1}}}{F^*{E_i}}$$
is the HN filtration of $F^*({E_{i+1}}/{E_i})$. 
Now it is enough to prove the 

\noindent{\bf Claim}.~  
$$0\subset E_{01}\subset \cdots \subset E_{0t_0} \subset F^*{
E_1}\subset \cdots \subset F^*{E_i} \subset E_{i1} \subset \cdots
\subset E_{it_i} \subset F^*{E_{i+1}} \subset \cdots \subset F^*V
$$
is the HN filtration of $F^*V$. 

\noindent~{\underline{Proof of the claim}}.~By construction, for  $0\leq
i\leq l$ and for  $1\leq j< t_i$, we have

$$\mu\left(\frac{E_{ij}}{E_{i(j-1)}}\right) >
\mu\left(\frac{E_{i(j+1)}}{E_{ij}}\right)
$$
and 
$$\frac{E_{ij}}{E_{i(j-1)}}, ~\frac{F^*{E_i}}{E_{i(t_i-1)}}
~~\mbox{and}~~  
\frac{E_{i1}}{F^*{E_i}}$$
 are semistable. 
Hence, by Definition~\ref{d2}, it is enough to prove that
$$\mu\left(\frac{F^*{E_i}}{E_{i-1t_{i-1}}}\right) >
\mu\left(\frac{E_{i1}}{F^*{E_i}}\right),$$
Now, by Corollary~$2^p$ of 
[SB],
we have 
\begin{equation}\label{e1}
0\leq \mu_{max}F^*\left(\frac{{E_{i+1}}}{{E_i}}\right)
-\mu_{min}F^*\left(\frac{{E_{i+1}}}{{E_i}}\right) \leq
(2g-2)(r-1).
\end{equation} 
By Remark~\ref{r3}, for all $0\leq i\leq l$, we have 
$$\mu_{max}F^*\left(\frac{{E_{i+1}}}{{E_i}}\right) \geq
\mu(F^*\left(\frac{{E_{i+1}}}{{E_i}}\right)) \geq 
\mu_{min}F^*\left(\frac{{E_{i+1}}}{{ E_i}}\right).$$
Therefore 
$$0\leq \mu_{max}F^*\left(\frac{{E_{i+1}}}{{E_i}}\right)-
\mu(F^*\left(\frac{{E_{i+1}}}{{E_i}}\right)) \leq
(2g-2)(r-1),$$
which means
\begin{equation}\label{e2}
0\leq \mu\left(\frac{E_{i1}}{F^*{E_i}}\right) -
p\mu_{i+1}\leq (2g-2)(r-1).
\end{equation}
Similarly
$$0\leq \mu(F^*\left(\frac{{ E_i}}{{E_{i-1}}}\right))-
\mu_{min}(F^*\left(\frac{{ E_{i}}}{{ E_{i-1}}}\right))
\leq (2g-2)(r-1)$$
which means
\begin{equation}\label{e3}
0\leq 
p\mu_i - \mu\left(\frac{F^*{
E_i}}{E_{i-1t_{i-1}}}\right) \leq (2g-2)(r-1). \end{equation}

Now, multiplying (\ref{e2}) and (\ref{e3}) by $-1$ and adding,
we
get  
\begin{equation}\label{eqq}
-4(g-1)(r-1)+p(\mu_i-\mu_{i+1})\leq  \mu\left(\frac{F^*{
E_i}}{E_{i-1t_{i-1}}}\right) - \mu\left(\frac{E_{i1}}{F^*{
E_i}}\right) \leq
p(\mu_i-\mu_{i+1}).
\end{equation}
Since  $p >
4(g-1)r^3$, Lemma~\ref{l0} implies that
$$-4(g-1)(r-1)+p(\mu_i-\mu_{i+1}) > 0,$$
and hence 
$$\mu\left(\frac{F^*{E_i}}{E_{i-1t_{i-1}}}\right) >
\mu\left(\frac{E_{i1}}{F^*{ E_i}}\right),$$
This proves the claim, and hence the lemma.\end{proof}
 
\begin{defn}\label{d5}A filtration by subbundles
\[0=E_0\subset E_1\subset\cdots\subset E_l\subset E_{l+1}=V\]
of $V$ is a {\em strongly semistable HN
filtration} if
\begin{enumerate}
\item it is the HN filtration and
\item $E_1, E_2/E_1, \ldots, E_{l+1}/E_l$ are strongly semistable vector
bundles. 
\end{enumerate} \end{defn} 

\begin{remark}\label{r0}If the HN filtration  
\[0=E_0\subset E_1\subset\cdots\subset E_l\subset E_{l+1}=V\]
 of $V$ is strongly semistable then, for any $k\geq 0$, the filtration
\[0=E_0\subset F^{k*}E_1\subset\cdots\subset F^{k*}E_l\subset
F^{k*}E_{l+1}=F^{k*}V\]
is the strongly semistable HN filtration of  $F^{k*}V$.\end{remark}

\begin{remark}\label{r2}Note that, if $\rank~V =r$ and ${\rm  char}~k =
p>4(g-1)r^3$, 
then Lemma~\ref{l1} implies that there exists $s\geq 0$ such that the HN
filtration of $F^{s*}V$ is strongly semistable. Therefore, Theorem~4.5
of [L] follows in this case.
\end{remark}

\begin{defn}\label{d4} Let $E$ be a vector bundle on $X$. A 
vector bundle $F_j\neq 0$ occuring in the HN filtration of $F^{s*}E$ is
said
to {\it almost descend} to a  bundle $E_{i1}$ occuring in the HN
filtration of  
$E$ if $F_j\subseteq F^{s*}E_{i1}$ and  $E_{i1}$ is the
smallest bundle in the
HN filtration of $E$, with this property.
\end{defn}

\begin{remark}\label{r8}Note that, if $p >4(g-1)({\rank}~E)^3$, then
by Lemma~\ref{l1}, we have the following transitivity property: if $k
\leq
s$ such
that 
$F_j$
almost descends to a bundle ${\tilde E}_{i1}$ in the HN filtration of
$F^{(s-k)*}E$, and ${\tilde E}_{i1}$ almost descends to a bundle
$E_{ik}$ occuring in the HN filtration of $F^{k*}E$, then
 $F_j$ almost descends to the  bundle
$E_{ik}$. \end{remark}

\begin{lemma}\label{l5}Let $E$ be a vector bundle on $X$ of rank $r$
and let the characteristic $p$ satisfy
$p>4(g-1)r^3$. Let 
$F_j\neq 0$ be a subbundle in the HN filtration of $F^{s*}E$, which
almost descends to a vector bundle $E_i$ occuring in the HN filtration
of 
$E$. Then 
$$ \frac{{\underline
\mu}(F_j)}{p^s} = {\underline \mu}(E_i) +\frac{C}{p},$$
where $|C|\leq 4(g-1)(r-1)$ and ${\underline \mu(F_j)}$ is as defined in
Notation~\ref{n1}.\end{lemma}
\begin{proof}Let $F_{j-1}$ be the vector bundle on $X$ such that
$F_{j-1}\subset F_j$ are two consecutive subbundles of the HN filtration
of
$F^{s*}E$. Therefore, by Lemma~\ref{l1}, there exist two consecutive
subbundles 
$E_{i_1-1}\subset E_{i_1}$ in the HN filtration of $F^{(s-1)*}E$ 
 such that 
$$F^*E_{i_1-1} \subseteq F_{j-1} \subset F_j\subseteq F^*E_{i_1}.$$
In
particular, we are in the situation that $E_{i_1}/E_{i_1-1}$ is
semistable on $X$ and 
\begin{enumerate}
\item either $F_{j-1}/F^*E_{i_1-1} = 0$ in $F^*(E_{i_1}/E_{i_1-1})$,
and $F_j/F^*E_{i_1-1}$ is the first nonzero vector bundle in the HN
filtration of $F^*(E_{i_1}/E_{i_1-1})$ or
\item $F_{j-1}/F^*E_{i_1-1}\subset F_j/F^*E_{i_1-1}$ are
two consecutive subundles in the HN filtration of
$F^*(E_{i_1}/E_{i_1-1})$.
\end{enumerate}
In both the cases, by Definition~\ref{d2}, we
have 
$$\mu_{\rm min}F^*\left(\frac{E_{i_1}}{E_{i_1-1}}\right)\leq
\mu\left(\frac{F_j}{F_{j-1}}\right) \leq 
 \mu_{\rm max}F^*\left(\frac{E_{i_1}}{E_{i_1-1}}\right).$$
Therefore, Corollary~$2^p$ of [SB] implies 
$$-2(g-1)(r-1)\leq {\underline \mu}(F_j)-
\mu\left(F^*\left(\frac{E_{i_1}}{E_{i_1-1}}\right)\right)\leq 
2(g-1)(r-1) $$
Note that 
$\mu\left(F^*\left(E_{i_1}/E_{i_1-1}\right)\right)
= p{\underline \mu}(E_{i1})$. Therefore we have 
${\underline \mu}(F_{j}) = p{\underline \mu}(E_{i1})+C_1,$
where $|C_1| \leq 2(g-1)(r-1)$.

Note $E_{i_1}$ is a nonzero subbundle in the HN filtration of
$F^{(s-1)*}E$ which  almost descends to $E_i$ occuring in the HN
filtration of
$E$. Hence, inductively one can prove that 
 $$ {\underline \mu}(E_{i1}) = p^{s-1}{\underline
\mu}(E)+p^{s-2}C_s+\cdots
+C_2,$$
where $|C_2|, 
\ldots, |C_s|
\leq 2(g-1)(r-1)$. 
Therefore 
$${\underline \mu}(F_j)  =  p^s{\underline
\mu}(E_i)+ p^{s-1}C_s+\cdots
+pC_2+C_1.$$
Therefore 
$$ \frac{{\underline \mu}(F_j)}{p^s} = {\underline \mu}(E_i)+
\frac{1}{p^s}(p^{s-1}C_s+\cdots
+pC_2+C_1).$$
But 
$$|(p^{s-1}C_s+\cdots
+pC_2+C_1)|
\leq  (1+\cdots+p^{s-1})(2(g-1)(r-1)).$$
Since $(1+p+\cdots+p^{s-1})/p^{s-1} \leq 2$, we have 
$$\frac{|p^{s-1}C_s+\cdots
+pC_2+C_1|}{p^{s-1}} \leq 4(g-1)(r-1).$$
Therefore we conclude that 
$$\frac{{\underline \mu}(F_j)}{p^s} = {\underline
\mu}(E_i) +\frac{C}{p},$$
where $|C|\leq 4(g-1)(r-1)$. This proves the lemma.\end{proof}

\begin{notation}\label{n3} Henceforth we assume that the characteristic
$p$ satisfies $p>4(g-1)r^3$.
We also fix  a vector bundle $V$ on $X$ of 
 rank $r$ with the HN filtration
$$0 ={E_0} \subset { E_1}\subset {E_2} \subset \cdots
\subset {
E_l}\subset { E_{l+1}} = V,$$
and let
$$\mu_i(V)= \mu\left(\frac{ E_i}{ E_{i-1}}\right) ~~~\mbox{and}~~~r_i(V)=
\rank\left(\frac{ E_i}{ E_{i-1}}\right).$$
Let
\begin{equation}\label{e7}
0\subset F_1\subset F_2\subset \cdots \subset F_t
\subset F_{t+1} =
F^{k*}V\end{equation}
be the HN filtration, and let 
$${\mu_i}(F^{k*}V)=\mu\left(\frac{F_i}{F_{i-1}}\right),~~
{r_i}(F^{k*}V) = \rank\left(\frac{F_i}{F_{i-1}}\right)
 ~~~\mbox{and}~~a_i(F^{k*}V) =
p^{-k}\mu\left(\frac{F_i}{F_{i-1}}\right).$$
\end{notation}

\begin{proposition}\label{8}With the notation as above, where
$p>4(g-1)r^3$,
if a vector bundle $F_j$ of the HN filtration of $F^{k*}V$
almost descends to a 
vector bundle ${E_i}$ of the HN filtration of $V$ then,
for any $m\geq 1$,
 $$a_j(F^{k*}V)^m = \mu_i(V)^m + \frac{C}{p},$$
where $|C|\leq 8gr({\mbox max}\{2|\mu_1|, \ldots, 2|\mu_{l}|,
2\}^{m-1})$.\end{proposition}
\begin{proof}By Lemma~\ref{l5}, we have
$$a_j(F^{k*}V) = \mu_{i}(V) + \frac{C_{ij}}{p},$$
where $|C_{ij}|\leq 4(g-1)(r-1)$.
Therefore 
$$a_j(F^{k*}V)^m-\mu_i(V)^m = {{m}\choose{1}}\mu_i(V)^{m-1}C_{ij}+\cdots
+{{m}\choose{m-1}}\mu_i(V)C_{ij}^{m-1}+{{m}\choose{m}}C_{ij}^m.$$
Now, as $|C_{ij}|/p\leq 1$, one can check that 
$$|a_j(F^{k*}V)^m-\mu_i(V)^m|\leq  
\frac{C}{p},$$
where $|C|\leq 8gr({\mbox max}\{2|\mu_1|, \ldots, 2|\mu_{l}|,
2\}^{m-1})$. This proves the proposition.\end{proof}

\section{Applications}

We extend Notation~\ref{n3} to the case, when the underlying field is of 
arbitrary characteristic,
as
follows.

\begin{notation}\label{n4}Let $X$ be a nonsingular curve over an
algebraically
closed field $k$ and $V$ a vector bundle on $X$, with HN filtration
$$0 = E_0 \subset E_1 \subset \ldots \subset E_l\subset E_{l+1} = V.$$ 
\begin{enumerate}
\item If ${\rm char}~k =p>0$, then we define the numbers $\mu_i(V)$,
$r_i(V)$,
${\mu}_i(F^{k*}V)$, ${r_i}(F^{k*}V)$ and ${a_i}(F^{k*}V)$ as in
Notation~\ref{n3}. 
Moreover, we 
 choose  an integer $s\geq 0$  such that
$F^{s*}(V)$ has a strongly semistable HN filtration and we denote 
$${\tilde a_i}(V) = a_i(F^{s*}(V))~~\mbox{and}~~ {\tilde r_i}(V) =
r_i(F^{s*}(V))$$
(note that, by Remark~\ref{r0}, 
these numbers are 
 independent of the choice of such an $s$),   
\item If ${\rm char}~k = 0$, define 
${\tilde a_i}(V) = \mu_i(V) = {\mu}(E_i/E_{i-1})$, and ${\tilde
r_i}(V) = r_i(V) = \rank(E_i/E_{i-1})$.
\end{enumerate}
\end{notation}

\begin{proposition}\label{p1} Let $X_A\by{} {\rm Spec}~A$ be 
a projective morphism of Noetherian schemes, smooth of relative
dimension~1, where $A$ is a finitely
generated $\Z$-algebra and is an integral domain. Let $\sO_{X_A}(1)$ be an 
$f$-very ample invertible sheaf on $X_A$. Let $V_A$ be a vector
bundle on $X_A$. For  $s\in {\rm Spec}~A$, let $V_s= V_A\tensor_A
{\overline {k(s)}}$ be the induced vector bundle on the smooth projective 
curve $X_s = X_A\tensor_A{\overline{k(s)}}$.
Let $s_0 = {\rm Spec}~Q(A)$ be the generic point of  ${\rm
Spec}~A$.
Then, 
\begin{enumerate}
\item for any $k\geq 0$ and $m\geq 0$, we have 
$$\lim_{s\to s_0}\sum_j{r_j}(F^{k*}V_s) a_j(F^{k*}V_s)^m =
\sum_i{r_i}(V_{s_0})\mu_i(V_{s_0})^m.$$
\item Similarly
$$\lim_{s\to s_0}\sum_j{\tilde r_j}(V_s) {\tilde a_j}(V_s)^m =
 \sum_i{r_i}(V_{s_0})\mu_i(V_{s_0})^m,$$
\end{enumerate}
where in both the limits, $s$ runs over closed points of ${\rm{Spec}}~A$. 
\end{proposition} 
\begin{proof} To prove the proposition, one can replace ${\rm Spec}~A$ by
an affine open subset. 
We may  assume that
 the HN filtration $\{E_i\}$ of $V_{s_0}$ on $X_{s_0}$ is defined
on the model
 $V_A$, and restricts to a filtration of 
$V_s$ by submodules, for each $s$. 
Under this reduction, the slopes of the respective quotients are preserved. 
Finally, by an openness property of semistable vector bundles (\cite{Ma}),
we may assume (after localizing $A$ if necessary) that the resulting 
filtration of $V_s$ on $X_s$ is the HN filtration of $V_s$.  
Therefore we can choose $A$ such that, for any closed point $s\in {\rm
Spec}~A$, we have
$${\rm {char}}~k(s) > 4(\mbox{genus}~X_s-1)(\rank~V_s)^3 = 
4(\mbox{genus}~X_{s_0}-1)(\rank~V_{s_0})^3.$$
Therefore, if denote
$$M =  8({\rm genus}(X_{s_0}))r(V_{s_0})(\max\{2,
2|\mu_1(V_{{s_0}})|,\ldots,|\mu_l(V_{s_0})|\}^{m-1}),$$
 then, by
Proposition~\ref{8},
we have 
$$\begin{array}{lcl}
\sum_j{r}_j(F^{k*}V_{X_s})a_j(F^{k*}V_{X_s})^m
& = & \sum_i{r}_i(V_{X_s})\left(\mu_i(V_{X_s})^m
+\displaystyle{\frac{C_i}{p}}\right),~~~
\mbox{where}~~|C_i|\leq M\\
& = & \sum_i{r}_i(V_{X_{s_0}})\mu_i(V_{X_{s_0}})^m
+\displaystyle{\frac{C_{s_k}}{p}},\end{array}
$$
where 
$|C_{s_k}|\leq r(V_{s_0})M$. In particular, for every closed point $s\in
{\mbox{Spec}}~A$, we have 
$$\sum_j{\tilde r_j}(V_{X_s}){\tilde a_j}(V_{X_s})^m
= \sum_i{r}_i(V_{X_{s_0}})\mu_i(V_{X_{s_0}})^m
+\displaystyle{\frac{C_s}{p}},$$
where $|C_s|\leq r(V_{s_0})M$. 
Now the proposition follows easily.
\end{proof}

\begin{cor}\label{c}Along with Notation~\ref{n4}, if we 
 denote (as defined in [B2]), for
 ${\rm char}~k >0$, $\mu_{HK}(V) = \sum_i{\tilde r_i}(V)a_i(V)^2$, and 
for  ${\rm char}~k = 0$, $\mu_{HK}(V) = \sum_j{r_j}(V)\mu_j(V)^2$,
 then 
$$\lim_{s\to s_0}\mu_{HK}(V_{s})
= \mu_{HK}(V_{{s_0}}).$$
\end{cor}
\begin{proof}The corollary follows by substituting $m=2$ in the second
statement of
Proposition~\ref{p1}. \end{proof}

Let $k$ be an algebraically closed field of characteristic $0$. Let $R$ be
 a finitely generated $\N$-graded two dimensional domain over
$k$. Let $I\subset R$ be an
homogeneous ideal of finite colength. Then there exists a finitely
generated $\Z$-algebra $A\subseteq k$, a finitely generated $\N$-graded
algebra $R_A$
over $A$ and an homogeneous ideal $I_A\subset R_A$ such that 
$R_A\tensor_Ak = R$ and for any closed point $s\in {\rm Spec}~A$ ({\it
i.e.}
maximal ideal of $A$) the ring $R_s = R_A\tensor_A k(s)$ is a finitely
generated $\N$-graded 2-dimensional
 domain (which is a normal domain if $R$ is normal) over
$k(s)$ and the ideal $I_s = {\rm{Im}}(I_A\tensor_A k(s)) \subset R_s $ is
an
homogeneous ideal of finite colength.

Moreover, if, for the pair $(R, I)$,  we have a {\it spread} $(A,R_A,
I_A)$
as above and
 $A\subset A' \subset
k$, for some finitely generated $\Z$-algebra $A'$  then $(A', R_{A'},
I_{A'})$
satisfy the same properties as $(A, R_A, I_A)$. Hence we may always assume
that the spread $(A, R_A, I_A)$ as above is chosen 
such that $A$ contains a given finitely generated algebra $A_0 \subseteq
k$. 

\begin{thm}\label{t2}Let $R$ be a standard graded two dimensional
domain over $k$. Let $I \subset R$ be an homogeneous ideal of finite
colength. Let $(A, R_A, I_A)$ be a spread as given above. Then 
$$\lim_{s\to s_0}HKM(R_{s}, I_s)$$
 exists,
where $s_0 = {\rm Spec}~Q(A)$ is the generic point of  ${\rm
Spec}~A$, and the limit is taken over closed points $s\in {\rm{Spec}}A$. 
\end{thm} 
\begin{proof}Let $R\by{}S$ be the normalization of $R$. Then $R\by{} S$
is a finite graded map of degree $0$, and $Q(R) = Q(S)$, such
that $S$ is a
finitely generated $\N$-graded 2-dimensional normal domain over $k$.
Now, for pairs $(R, I)$, $(S, IS)$, we choose spreads $(A,R_A, I_A)$
and $(A, S_A , IS_A)$ such that for every closed point $s\in
{\mbox{Spec}}~A$,
the
natural map $R_s = R_A\tensor k(s)\by{} S_s =S_A\tensor k(s)$ is
a finite graded map of degree $0$. Therefore we have the following
commutative
diagrams of finite
horizontal maps 

$$\begin{array}{ccc}
{\rm Proj}~R & \longleftarrow & {\rm Proj}~S\\
\downarrow & & \downarrow\\
{\rm Proj}~R_A & \longleftarrow & {\rm Proj}~S_A.\end{array}$$
It follows that, for every  $s\in {\rm Spec}~A$, the
corresponding map of curves 
$${\rm Proj}~R_A\tensor_A k(s)\to {\rm Proj}~S_A\tensor_Ak(s)$$ is a
finite
map, where the curve ${\rm Proj}~S_A\tensor_A k(s)$ is
nonsingular. 
Therefore 
$$HKM(R_s, I_s) = HKM(S_s, IS_s), ~~\mbox{for every closed point}~~{s\in
{\rm 
{Spec}}~A}.$$
Therefore it is enough to prove the following 

\vspace{5pt}

\noindent{\bf Claim}. $\lim_{s\to s_0}HKM(S_s, IS_s)$ exists.

\noindent{\underline{Proof of the claim}}: Let  $I$ and $IS_A$  be
generated by the
set $\{f_1,\ldots, f_k\}$, where $\deg~f_i = d_i$.
 We
have a short exact sequence of $\sO_{X_A}$-sheaves (see [B1] and [T1]):
\begin{equation}\label{e5}
0\by{} V_{A}\by{}
\oplus_{i=1}^k\sO_{X_A}(1-d_i)\by{}
\sO_{X_A}(1)\by{} 0
\end{equation}
where $\sO_{X_A}(1-d_i)\by{}\sO_{X_A}(1)$ is multiplication by
$f_i$. 
Restricting (\ref{e5}) to the fiber $X_s$, we get
$$0\by{} V_{s}=V_{A}\tensor_Ak(s) \by{}
\oplus_{i=1}^k\sO_{X_s}(1-d_i)\by{} \sO_{X_s}(1)\by{} 0.$$ 
Note that (see [B1] and [T1]),
$$HKM(S_s, IS_s) = 
\frac{\deg~{\rm Proj}~S_s}{2}\left(\sum_i{\tilde
r_i}(V_s)a_i(V_s)^2)-\sum_{i=1}^{k}d_i^2\right) $$
Therefore 
$$\lim_{s\to s_0}HKM(S_s, IS_s) = 
\frac{\deg~{\rm Proj}~S}{2}\left(\lim_{s\to s_0}\sum_i{\tilde
r_i}(V_s)a_i(V_s)^2- \sum_{i=1}^{k}d_i^2\right),$$
which,  by Proposition~\ref{p1}, exists.
This proves the theorem.
\end{proof}

\begin{remark}\label{r6} Let $R$ be a standard graded $2$ dimensional
domain
over a field of characteristic $0$. Let $I \subset R$ be a homogeneous
ideal
of finite colength. Then for the pair $(R,I)$ we choose a spread $(A, X_A,
I_A)$ as described earlier and define 
\begin{equation}\label{e9}
HKM(R,I) = \lim_{s\to s_o}HKM(R_s, I_s).
\end{equation}
This is,  inherently, a well defined notion (${\it i.e.}$, irrespective of
a
choice of
generators of $I$), since in positive characteristic $HKM(R_s, I_s)$ is
independent of a choice of generators of $I_s$.
We extend this definition  to a standard graded $2$-dimensional 
ring $R$, over a field $k$ of  characteristic $0$, and a homogeneous ideal
$I\subset R$ of finite colength as 
$$HKM(R, I) = \sum _{{\bf p}\in{\rm Spec}R, \dim R/{\bf p}
=2}\ell_{R_{\bf p}}(R_{\bf p})HKM(R/{\bf p},
IR/{\bf p}),$$
Note that a notion of $HKM(R, I)$, when $R$ is also a normal domain (${\it
i.e.}$,
${\rm Proj}~R$ is a smooth curve) over a field of characteristic $0$, is
given in
[B2]
as 
\begin{equation}\label{e8}
HKM(R, I) = \frac{\deg~{\rm Proj}~R}{2}\left(\mu_{HK}(V)-
\sum_{i=1}^kd_i^2\right),
\end{equation}
where $V$ is the vector bundle given by 
$$0\by{} V\by{} \oplus_{i} \sO_X(1-d_i)\by{} \sO_{X}(1)\by{} 0.$$
By Corollary~\ref{c}, these two definitions
(\ref{e8}) and (\ref{e9}) coincide, in this case.  
\end{remark}

\begin{remark}\label{r7}It follows from Remark~4.13 of [T1] that, for
every closed point $s$ in ${\mbox{Spec}}~A$, where $(A, R_A, I_A)$ is a
spread for the
pair $(R, I)$, we have 
$$HKM(R_s, I_s)\geq HKM(R, I),$$ and 
$HKM(R_s, I_s) = HKM(R, I)$ if and
only if HN
filtration of $V_s$ is the strongly semistable HN filtration. Note that  
$HKM(R, I)$ is a rational number, expressed  in terms of  the slopes
of the subquotients of HN filtration of $V_s$. For example, when $V_s$ is
semistable then
$$HKM(R, I) = \deg(R_s)/2\left((\sum_id_i)^2/(t-1)-\sum_id_i^2\right).$$
\end{remark}
\begin{remark}\label{r9}As observed in the above remark, 
$$\{HKM(R_s, I_s)- HKM(R,I)\mid s \in \{\mbox{closed points of}~
{\rm{Spec}}A\}\}$$ is a sequence of positive rational numbers (indexed by
the closed points), 
converging to $0$. The   
following example of Monsky implies that  it could be {\it oscillating}.

Let 
$$R_p = k[X,Y,Z]/(X^4+Y^3Z+Z^3X), ~\mbox{where}~{\rm{char}.}k = p.$$ Then 
$$\begin{array}{lcl}
 HKM(R_p,(X,Y,Z)R_p) & = &
3+\displaystyle{\frac{1}{4p^2}},~\mbox{if}~p\equiv
\pm4(9)\\
  & = & 3+\displaystyle{\frac{1}{4p^4}},~\mbox{if}~p\equiv \pm2(9)\\
 & = & 3,~\mbox{if}~p\equiv \pm 1(9).
\end{array}$$
Now, let $X_p= {\rm{Proj}}~R_p$.
If we consider the short exact sequence
$$ 0\by{} V_{X_p}\by{} \sO_{X_p}\oplus\sO_{X_p}\oplus\sO_{X_p}
\by{}\sO_{X_p}(1)\by{} 0,$$
 where the second map is given by
$(f_1, f_2, f_3) \to Xf_1+Yf_2+Zf_3$.

We also recall the following result of [T2] 
\begin{cor}\label{20}Let $X_p= {\mbox{Proj}}~(R_p = k[X,Y,Z]/(f))$ be
a nonsingular plane
curve of degree $d$ over an algebraically closed field $k$ of
characteristic
$p>0$. 
Then $$ HKM(X_p,\sO_{X_p}(1)) = HKM(R_p, (X,Y,Z)R_p) =
\frac{3d}{4}+\frac{l^2}{4dp^{2s}},$$
where $s\geq 1$ is a number such that $F^{(s-1)*}V_{X_p}$ is
semistable and
$F^{s*}V_{X_p}$ is not semistable (if $F^{t*}V_{X_p}$ is
semistable for all $t\geq 0$, we take
$s=\infty$) and $l$ is an integer congruent to $pd$ (mod 2) with
$0\leq l \leq d(d-3)$.\end{cor}

This implies that in the example given above
\begin{enumerate}
\item if $p\equiv \pm-4(9)$ and $p>>0$ then 
$l=2$ and $s=1$, {\it i.e.}
$V_{X_p}$ is semistable, and $F^{*}(V_{X_p})$ is not semistable and has
strongly semistable HN
filtration and
$$a_1(V_{X_p}) = \mu(V_{X_p}) + \frac{1}{p} ~~\mbox{and}~~
a_2(V_{X_p}) = \mu(V_{X_p}) - \frac{1}{p}$$
In particular $\mu_{HK}(V_{X_p}) = 2\mu(V_{X_p})^2 +\frac{2}{p^2}$.

\item if $p\equiv \pm-2(9)$ and $p>>0$ then 
$l=2$ and $s=2$, {\it i.e.}
$F^*V_{X_p}$ is semistable, and $F^{2*}(V_{X_p})$ is not semistable and
has strongly semistable
HN
filtration and
$$a_1(V_{X_p}) = \mu(V_{X_p}) + \frac{1}{p^2} ~~\mbox{and}~~
a_2(V_{X_p}) = \mu(V_{X_p}) - \frac{1}{p^2}$$
In particular $\mu_{HK}(V_{X_p}) = 2\mu(V_{X_p})^2 +\frac{2}{p^4}$.
\end{enumerate}

In particular, for  $p>>0$ the numbers $a_1(V_{X_p})$
$a_2(V_{X_p})$ do not
eventually become constant or
a well defined function of $p$, but keep oscillating and converge to
$\mu(V_X)$.  
\end{remark}

\end{document}